\documentclass[11pt,a4paper]{article}
\usepackage{amssymb,latexsym,a4wide}
\def\exp{\mathop{\rm exp}\nolimits}

\def\proof{{\bf Proof}.\ }

\newtheorem{formula}{}[section]
\newtheorem{proposition}[formula]{Proposition}
\newtheorem{definition}[formula]{Definition}
\newtheorem{corollary}[formula]{Corollary}
\newtheorem{remark}[formula]{Remark}
\newtheorem{lemma}[formula]{Lemma}
\newtheorem{theorem}[formula]{Theorem}

\def\prop{\begin{proposition}}
\def\propl#1{\begin{proposition}\label{#1}}
\def\eprop{\end{proposition}}
\def\exp{\begin{example}}
\def\expl#1{\begin{example}\label{#1}}
\def\eexp{\end{example}}
\def\thrm{\begin{theorem}}
\def\thrml#1{\begin{theorem}\label{#1}}
\def\ethrm{\end{theorem}}
\def\rmrk{\begin{remark}}
\def\rmrkl#1{\begin{remark}\label{#1}}
\def\ermrk{\end{remark}}
\def\dfntn{\begin{definition}}
\def\dfntnl#1{\begin{definition}\label{#1}}
\def\edfntn{\end{definition}}
\def\nmrt{\begin{enumerate}}
\def\enmrt{\end{enumerate}}

\def\qtn{\begin{equation}}
\def\qtnl#1{\begin{equation}\label{#1}}
\def\eqtn{\end{equation}}
\def\lmm{\begin{lemma}}
\def\lmml#1{\begin{lemma}\label{#1}}
\def\elmm{\end{lemma}}
\def\crllr{\begin{corollary}}
\def\crllrl#1{\begin{corollary}\label{#1}}
\def\ecrllr{\end{corollary}}
\begin{document}
\title{A note on the computation of Puiseux series solutions of the Riccatti equation associated with a homogeneous linear ordinary differential equation}
\author{
Ali AYAD \\[-1pt]
\small IRMAR , Campus universitaire de Beaulieu \\[-3pt]
\small Universit\'e Rennes 1, 35042, Rennes, France\\[-3pt]
{\tt \small ali.ayad@univ-rennes1.fr}\\[-3pt]
\small And \\[-3pt]
\small IRISA (Institut de recherche en informatique et syst\`emes al\'eatoires)\\[-3pt]
\small INRIA-Rennes, Campus universitaire de Beaulieu \\[-3pt]
\small Avenue du G\'en\'eral Leclerc, 35042 Rennes Cedex, France\\[-3pt]
\small ali.ayad@irisa.fr \\[-3pt]
}
\date{}

\maketitle

\begin{abstract}
We present in this paper a detailed note on the computation of Puiseux series solutions of the Riccatti equation associated with a homogeneous linear ordinary differential equation. This paper is a continuation of~\cite{Aya7} which was on the complexity of solving arbitrary ordinary polynomial differential equations in terms of Puiseux series.
\end{abstract}

\section*{Introduction}

Let $K=\mathbb{Q}(T_1, \dots ,T_l)[\eta]$ be a finite extension of a finitely generated field over $\mathbb{Q}$. The variables $T_1, \dots ,T_l$ are algebraically independent over $\mathbb{Q}$ and $\eta$ is an algebraic element over the field $\mathbb{Q}(T_1, \dots ,T_l)$ with the minimal polynomial $\phi \in \mathbb{Z}[T_1, \dots ,T_l][Z]$. Let $\overline{K}$ be an algebraic closure of $K$ and consider the two fields:
$$L= \cup_{\nu \in \mathbb{N}^*}K((x^{\frac{1}{\nu}})), \quad \mathcal{L}= \cup_{\nu \in \mathbb{N}^*}\overline{K}((x^{\frac{1}{\nu}}))$$ 
which are the fields of fraction-power series of $x$ over $K$ (respectively $\overline{K}$), i.e., the fields of Puiseux series of $x$ with coefficients in $K$ (respectively $\overline{K}$). Each element $\psi \in L$ (respectively $\psi \in \mathcal{L}$) can be represented in the form $\psi =\sum_{i\in \mathbb{Q}}c_ix^i, \quad c_i\in K$ (respectively $c_i\in \overline{K}$). The order of $\psi$ is defined by $ord (\psi):= \min \{i\in \mathbb{Q}, c_i\neq 0\}$. The fields $L$ and $\mathcal{L}$ are differential fields with the differentiation  
$$\frac{d}{dx}(\psi)=\sum_{i\in \mathbb{Q}}ic_ix^{i-1}.$$
\newline
Let $S(y)=0$ be a homogeneous linear ordinary differential equation which is writen in the form 
$$S(y)=f_ny^{(n)}+\cdots +f_1y'+f_0y$$
where $f_i\in K[x]$ for all $0\leq i\leq n$ and $f_n\neq 0$ (we say that the order of $S(y)=0$ is $n$). Let $y_0, \dots , y_n$ be new variables algebraically independent over $K(x)$. We will associate to $S(y)=0$ a non-linear differential polynomial $R\in K[x][y_0, \dots ,y_n]$ such that $y$ is a solution of $S(y)=0$ if and only if $\frac{y'}{y}$ is a solution of $R(y)=0$ where the last equation is the ordinary differential equation $R(y, \frac{dy}{dx}, \dots , \frac{d^ny}{dx})=0$. We consider the change of variable $z=\frac{y'}{y}$, i.e., $y'=zy$, we compute the successive derivatives of $y$ and we make them in the equation $S(y)=0$ to get a non-linear differential equation $R(z)=0$ which satisfies the above property. $R$ is called the {\it Riccatti} differential polynomial associated with $S(y)=0$. We will describe all the fundamental solutions (see e.g. \cite{Was, Gri9}) of the differential equation $R(y)=0$ in $\mathcal{L}$ by a differential version of the Newton polygon process. There is another way to formulate $R$: let $(r_i)_{i\geq 0}$ be the sequence of the following differential polynomials
$$r_0=1, \, \, r_1=y_0, \dots , r_{i+1}=y_0r_i+Dr_i, \, \, \forall i\geq 1,$$
where $Dy_i=y_{i+1}$ for any $0\leq i\leq n-1$. We remark that for all $i\geq 1$, $r_i\in \mathbb{Z}[y_0, \dots ,y_{i-1}]$ has total degree equal to $i$ w.r.t. $y_0, \dots ,y_{i-1}$ and the only term of $r_i$ of degree $i$ is $y_0^i$.

\lmml{comb}
The non-linear differential polynomial 
$$R=f_nr_n+\cdots +f_1r_1+f_0r_0\in K[x][y_0, \dots ,y_n]$$
is the {\it Riccatti} differential polynomial associated with $S(y)=0$.
\elmm

\section{Newton polygon of $R$}

Let $R$ be the {\it Riccatti} differential polynomial associated with $S(y)=0$ as in Lemma~\ref{comb}. We will describe the Newton polygon $\mathcal{N}(R)$ of $R$ in the neighborhood of $x=+\infty$ which is defined explicitly in~\cite{Aya7}. For every $0\leq i\leq n$, mark the point $(\deg (f_i), i)$ in the plane $\mathbb{R}^2$. Let $\mathcal{N}$ be the convex hull of these points with the point $(-\infty , 0)$.

\lmml{tritro}
The Newton polygon of $R$ in the neighborhood of $x=+\infty$ is $\mathcal{N}$, i.e., $\mathcal{N}(R)=\mathcal{N}$.
\elmm
\proof For all $0\leq i\leq n$, $\deg_{y_0, \dots , y_{i-1}}(r_i)=i$ and the only term of $r_i$ of degree $i$ is $y_0^i$, then $lc (f_i)x^{\deg (f_i)}y_0^i$ is a term of $R$ and $\mathcal{N}\subset \mathcal{N}(R)$. For any other term of $f_ir_i$ in the form $bx^jy_0^{\alpha_0}\cdots y_{i-1}^{\alpha_{i-1}}$, where $b\in K, \, \, j<\deg (f_i)$ and $\alpha_0+\cdots +\alpha_{i-1}<i$, the corresponding point $(j-\alpha_1-\cdots -(i-1)\alpha_{i-1}, \alpha_0+\cdots +\alpha_{i-1})$ is in the interior of $\mathcal{N}$. Thus $\mathcal{N}(R)\subset \mathcal{N}$. $\Box$

\lmml{}
For any edge $e$ of $\mathcal{N}(R)$, the {\it characteristic} polynomial of $R$ associated with $e$ is a non-zero polynomial. For any vertex $p$ of $\mathcal{N}(R)$, the {\it indicial} polynomial of $R$ associated with $p$ is a non-zero constant. Moreover, if the ordinate of $p$ is $i_0$, then $h_{(R, p)}(\mu)=lc (f_{i_0})\neq 0$. 
\elmm
\proof
By Lemma~\ref{tritro} each edge $e\in E(R)$ joints two vertices $(\deg (f_{i_1}), i_1)$ and $(\deg (f_{i_2}), i_2)$ of $\mathcal{N}(R)$. Moreover, the set $N(R, a(e), b(e))$ contains these two points. Then
$$0\neq H_{(R, e)}(C)=lc (f_{i_1})C^{i_1}+lc (f_{i_2})C^{i_2}+ t.$$
where $t$ is a sum of terms of degree different from $i_1$ and $i_2$. For any vertex $p\in V(R)$ of ordinate $i_0$, $lc (f_{i_0})x^{\deg (f_{i_0})}y_0^{i_0}$ is the only term of $R$ whose corresponding point $p$. Then 
$$h_{(R, p)}(\mu)=lc (f_{i_0})\neq 0. \Box$$

\crllrl{chlala}
For any edge $e\in E(R)$, the set $A_{(R, e)}$ is a finite set. For any vertex $p\in V(R)$, we have $A_{(R, p)}=\emptyset$.
\ecrllr

\section{Derivatives of the Riccatti equation}

For each $i\geq 0$ and $k\geq 0$, the $k$-th derivative of $r_i$ is the differential polynomial defined by 
$$r_i^{(0)}:= r_i, r_i^{(1)}:= r_i':= {\partial r_i \over \partial y_0} \, \, \textrm{and}\, \, r_i^{(k+1)}:= (r_i^{(k)})'={\partial^{k+1} r_i \over \partial y_0^{k+1}}.$$

\lmml{chelef}
For all $i\geq 1$, we have $r_i'=ir_{i-1}$. Thus for all $k\geq 0$, $r_i^{(k)}= (i)_kr_{i-k}$, where $(i)_0:=1$ and $(i)_k:= i(i-1)\cdots (i-k+1)$.
\elmm
\proof
We prove the first item by induction on $i$. For $i=1$, we have $r_1'=1=1.r_0$. Suppose that this property holds for a certain $i$ and prove it for $i+1$. Namely, 
\begin{eqnarray*}
r_{i+1}' & = & (y_0r_i+Dr_i)' \\
& = & y_0r_i'+r_i+Dr_i' \\
& = & iy_0r_{i-1}+r_i+D(ir_{i-1}) \\
& = & i(y_0r_{i-1}+Dr_{i-1})+r_i \\
& = & ir_i+r_i \\
& = & (i+1)r_i. 
\end{eqnarray*}
The second item is just a result of the first one (by induction on $k$). $\Box$
\dfntnl{}
Let $R$ be the {\it Riccatti} differential polynomial associated with $S(y)=0$. For each $k\geq 0$, the $k$-th derivative of $R$ is defined by
$$R^{(k)}:={\partial^k R \over \partial y_0^k}= \sum_{0\leq i\leq n}f_ir_i^{(k)}.$$
\edfntn

\lmml{7elweh}
For all $k\geq 0$, we have $$R^{(k)}= \sum_{0\leq i\leq n-k}(i+k)_kf_{i+k}r_i.$$ 
\elmm
\proof
For all $i<k$, we have $r_i^{(k)}=0$ because $\deg_{y_0}(r_i)=i$. Then by Lemma~\ref{chelef}, we get
\begin{eqnarray*}
R^{(k)} & = & \sum_{k\leq i\leq n}f_ir_i^{(k)} \\
& = & \sum_{k\leq i\leq n}f_i(i)_kr_{i-k} \\
& = &  \sum_{0\leq j\leq n-k}(j+k)_kf_{j+k}r_j
\end{eqnarray*}
with the change $j=i-k$. $\Box$

\crllrl{tranta}
The $k$-th derivative of $R$ is the {\it Riccatti} differential polynomial of the following linear ordinary differential equation of order $n-k$
$$S^{(k)}(y):=\sum_{0\leq i\leq n-k}(i+k)_kf_{i+k}y^{(i)}.$$
\ecrllr
\proof
By Lemmas~\ref{comb} and~\ref{7elweh}. $\Box$

\section{Newton polygon of the derivatives of $R$}

Let $0\leq k\leq n$ and $R^{(k)}$ be the $k$-th derivative of $R$. In this subsection, we will describe the Newton polygon of $R^{(k)}$. Recall that $R^{(k)}$ is the $k$-th partial derivative of $R$ w.r.t. $y_0$, then by the section 2 of~\cite{Aya7}, the Newton polygon of $R^{(k)}$ is the translation of that of $R$ defined by the point $(0, -k)$, i.e., $\mathcal{N}(R^{(k)})=\mathcal{N}(R)+\{(0, -k)\}$. The vertices of $\mathcal{N}(R^{(k)})$ are among the points $(\deg (f_{i+k}), i)$ for $0\leq i\leq n-k$ by Lemma~\ref{7elweh}. Then for each edge $e_k$ of $\mathcal{N}(R^{(k)})$ there are two possibilities: the first one is that $e_k$ is parallel to a certain edge $e$ of $\mathcal{N}(R)$, i.e., $e_k$ is the translation of $e$ by the point $\{(0, -k)\}$. The second possibility is that the upper vertex of $e_k$ is the translation of the upper vertex of a certain edge $e$ of $\mathcal{N}(R)$ and the lower vertex of $e_k$ is the translation of a certain point $(\deg (f_{i_0}), i_0)$ of $\mathcal{N}(R)$ which does not belong to $e$. In both possibilities, we say that the edge $e$ is associated with the edge $e_k$.

\lmml{deriveh}
Let $e_k\in E(R^{(k)})$ be parallel to an edge $e\in E(R)$. Then the {\it characteristic} polynomial of $R^{(k)}$ associated with $e_k$ is the $k$-th derivative of that of $R$ associated with $e$, i.e., 
$$H_{(R^{(k)}, e_k)}(C)=H_{(R, e)}^{(k)}(C).$$
\elmm
\proof
The edges $e_k$ and $e$ have the same inclination $\mu_e=\mu_{e_k}$ and $N(R^{(k)}, a(e_k), b(e_k))=N(R, a(e), b(e))+\{(0, -k)\}$. Then
 \begin{eqnarray*}
H_{(R^{(k)}, e_k)}(C) & = & \sum_{(\deg (f_{i+k}), i)\in N(R^{(k)}, a(e_k), b(e_k))} (i+k)_klc(f_{i+k})C^i \\
& = &  \sum_{(\deg (f_j), j)\in N(R, a(e), b(e))} (j)_klc(f_j)C^{j-k} \\
& = &  H_{(R, e)}^{(k)}(C). \Box 
\end{eqnarray*}

\crllrl{}
For any edge $e_k\in E(R^{(k)})$, the set $A_{(R^{(k)}, e_k)}$ is a finite set, i.e., $H_{(R^{(k)}, e_k)}(C)$ is a non-zero polynomial. For any vertex $p_k\in V(R^{(k)})$, we have $A_{(R^{(k)}, p_k)}=\emptyset$.
\ecrllr
\proof
By Corollaries~\ref{tranta} and~\ref{chlala}. $\Box$

\section{Newton polygon of evaluations of $R$}

Let $R$ be the {\it Riccatti} differential polynomial associated with $S(y)=0$. Let $0\leq c\in \overline{K}, \, \, \mu \in \mathbb{Q}$ and $R_1(y)= R(y+cx^{\mu})$. We will describe the Newton polygon of $R_1$ for different values of $c$ and $\mu$. 

\lmml{momop}
$R_1$ is the {\it Riccatti} differential polynomial of the following linear ordinary differential equation of order less or equal than $n$
$$S_1(y):= \sum_{0\leq i\leq n}\frac{1}{i!}R^{(i)}(cx^{\mu}) y^{(i)}.$$
\elmm
\proof
It is equivalent to prove the following analogy of Taylor formula: 
$$R_1=\sum_{0\leq i\leq n}\frac{1}{i!}R^{(i)}(cx^{\mu}) r_i$$ which is proved in Lemma 2.1 of~\cite{Gri9}. $\Box$
\newline
\newline
Then the vertices of $\mathcal{N}(R_1)$ are among the points $(\deg (R^{(i)}(cx^{\mu}), i)$ for $0\leq i\leq n$. Thus the Newton polygon of $R_1$ is given by (Lemma 2.2 of~\cite{Gri9}):

\lmml{mabe}
If $\mu$ is the inclination of an edge $e$ of $\mathcal{N}(R)$, then the edges of $\mathcal{N}(R_1)$ situated above $e$ are the same as in $\mathcal{N}(R)$. Moreover, if $c$ is a root of $H_{(R, e)}$ of multiplicity $m>1$ then $\mathcal{N}(R_1)$ contains an edge $e_1$ parallel to $e$ originating from the same upper vertex as $e$ where the ordinate of the lower vertex of $e_1$ equals to $m$. If $m=\deg H_{(R,e)}$, then $\mathcal{N}(R_1)$ contains an edge with inclination less than $\mu$ originating from the same upper vertex as $e$.
\elmm

\rmrkl{lol8}
If we evaluate $R$ on $cx^{\mu}$ we get
$$R(cx^{\mu}) = \sum_{0\leq i\leq n}f_i \times (c^ix^{i\mu}+t),$$
where $t$ is a sum of terms of degree strictly less than $i\mu$. Then 
$$lc (R(cx^{\mu}))=\sum_{i\in B}lc (f_i)c^i=\sum_{(\deg (f_i), i)\in e}lc (f_i)c^i=H_{(R, e)}(c),$$
where 
\begin{eqnarray*}
B & := & \{0\leq i\leq n;\, \, \deg (f_i)+i\mu=\max_{0\leq j\leq n}(\deg (f_j)+j\mu; \, \, f_j\neq 0)\} \\
& = & \{0\leq i\leq n;\, \, (\deg (f_i), i)\in e\, \, \textrm{and} \, \, f_i\neq 0\}.
\end{eqnarray*}
\ermrk

\lmml{quatrequatre}
Let $\mu$ be the inclination of an edge $e$ of $\mathcal{N}(R)$ and $c$ be a root of $H_{(R, e)}$ of multiplicity $m>1$. Then
$$H_{(R_1, e_1)}(C)=H_{(R, e)}(C+c)$$
where $e_1$ is the edge of $\mathcal{N}(R_1)$ given by Lemma~\ref{mabe}. In addition, if $e'$ is an edge of $\mathcal{N}(R_1)$ situated above $e$ (which is also an edge of $\mathcal{N}(R)$ by Lemma~\ref{mabe}) then $H_{(R_1, e)}(C)=H_{(R, e)}(C)$.
\elmm
\proof
We have
\begin{eqnarray*}
H_{(R, e)}(C+c) & = & \sum_{m\leq k\leq n}\frac{1}{k!}H_{(R, e)}^{(k)}(c)C^k \\
& = & \sum_{m\leq k\leq n}\frac{1}{k!}H_{(R^{(k)}, e)}(c)C^k \\
& = & \sum_{m\leq k\leq n}\frac{1}{k!}lc (R^{(k)}(cx^{\mu})) C^k \\
& = & H_{(R_1, e_1)}(C) 
\end{eqnarray*}
where the first equality is just the Taylor formula taking into account that $c$ is a root of $H_{(R, e)}$ of multiplicity $m>1$. The second equality holds by Lemma~\ref{deriveh}. The third one by Remark~\ref{lol8}. The fourth one by Lemma~\ref{momop} and by the definition of the {\it characteristic} polynomial. $\Box$

\section{Application of Newton-Puiseux algorithm to $R$}

We apply the Newton-Puiseux algorithm described in~\cite{Aya7} to the {\it Riccatti} differential polynomial $R$ associated with the linear ordinary differential equation $S(y)=0$. This algorithm constructs a tree $\mathcal{T}=\mathcal{T}(R)$ with a root $\tau_0$. For each node $\tau$ of $\mathcal{T}$, it computes a finite field $K_{\tau}=K[\theta_{\tau}]$, elements $c_{\tau}\in K_{\tau}\, \, \mu_{\tau}\in \mathbb{Q}\cup \{-\infty , +\infty\}$ and a differential polynomial $R_{\tau}$ as above. Let $\mathcal{U}$ be the set of all the vertices $\tau$ of $\mathcal{T}$ such that either $\deg (\tau)=+\infty$ and for the ancestor $\tau_1$ of $\tau$ it holds $\deg (\tau_1)<+\infty$ or $\deg (\tau)<+\infty$ and $\tau$ is a leaf of $\mathcal{T}$. There is a bijective correspondance between $\mathcal{U}$ and the set of the solutions of $R(y)=0$ in the differential field $\mathcal{L}$. The following lemma is a differential version of Lemma 2.1 of~\cite{Chi86} which separates any two different solutions in $\mathcal{L}$ of the {\it Riccatti} equation $R(y)=0$.

\lmml{separate}
Let $\psi_1, \psi_2\in \mathcal{L}$ be two different solutions of the differential {\it Riccatti} equation $R(y)=0$. Then there exist an integer $\gamma=\gamma_{12}$, $1\leq \gamma < n$, elements $\xi_1, \xi_2\in \overline{K}$, $\xi_1\neq \xi_2$ and a number $\mu_{12}\in \mathbb{Q}$ such that 
$$ord (R^{(\gamma)}(\psi_i)-\xi_ix^{\mu_{12}})<\mu_{12}, \, \, \textrm{for}\, \,  i=1, 2.$$
\elmm
\proof
By the above bijection, there are two elements $u_1$ and $u_2$ of $\mathcal{U}$ which correspond respectively to $\psi_1$ and $\psi_2$. Let $i_0=\max \{i\geq 0; \tau_i(u_1)=\tau_i(u_2)\}$. Denote by $\tau:=\tau_{i_0}(u_1)=\tau_{i_0}(u_2)$ and $\tau_1:=\tau_{i_0+1}(u_1)$, $\tau_2:=\tau_{i_0+1}(u_2)$. We have $\tau_1\neq \tau_2$ and $\epsilon :=\max (\mu_{\tau_1}, \, \, \mu_{\tau_2})$ is the inclination of a certain edge $e$ of $\mathcal{N}(R_{\tau})$. There are three possibilities for $\epsilon$:
\newline

- If $\mu_{\tau_2}<\mu_{\tau_1}$ then $\epsilon=\mu_{\tau_1}=\mu_e$. We have $c_{\tau_1}$ is a root of $H_{(R_{\tau}, e)}$ of multiplicity $m_1\geq 1$ and $R_{\tau_1}=R_{\tau}(y+c_{\tau_1}x^{\mu_{\tau_1}})$. Then by Lemma~\ref{mabe} there is an edge $e_1$ of $\mathcal{N}(R_{\tau_1})$ parallel to $e$ (so its inclination is $\epsilon=\mu_{\tau_1}$) originating from the same upper vertex as $e$ where the ordinate of the lower vertex of $e_1$ equals to $m_1$. In addition, $e$ is also an edge of $\mathcal{N}(R_{\tau_2})$ and by Lemma~\ref{quatrequatre}, we have $H_{(R_{\tau_2}, e)}(C)=H_{(R_{\tau}, e)}(C)$ and  
\qtnl{gho1}
H_{(R_{\tau_1}, e_1)}(C)=H_{(R_{\tau}, e)}(C+c_{\tau_1}).
\eqtn

- If $\mu_{\tau_1}<\mu_{\tau_2}$ then $\epsilon=\mu_{\tau_2}$. Then by Lemma~\ref{mabe} there is an edge $e_2$ of $\mathcal{N}(R_{\tau_2})$ parallel to $e$ originating from the same upper vertex as $e$. By Lemma~\ref{quatrequatre}, we have $H_{(R_{\tau_1}, e)}(C)=H_{(R_{\tau}, e)}(C)$ and 
\qtnl{gho2}
H_{(R_{\tau_2}, e_2)}(C)=H_{(R_{\tau}, e)}(C+c_{\tau_2}).
\eqtn

- If $\mu_{\tau_1}=\mu_{\tau_2}=\epsilon$ then $c_{\tau_1}$ and $c_{\tau_2}$ are two dictinct roots of the same polynomial $H_{(R_{\tau}, e)}(C)$. Then equalities of type~(\ref{gho1}) and~(\ref{gho2}) hold.
\newline
\newline
Set $\gamma =\deg_C (H_{(R_{\tau}, e)})-1\leq \deg_{y_0, \dots ,y_n}(R)-1\leq n-1<n$ and $\gamma \geq 1$ because that the polynomial $H_{(R_{\tau}, e)}(C)$ has at least two distinct roots $c_{\tau_1}$ and $c_{\tau_2}$. Moreover, we have $ord (\psi_i-y_{\tau_i})<\epsilon$ for $i=1, 2$. Let $\overline{\xi}_1\in K_{\tau_1}$ and $\overline{\xi}_2\in K_{\tau_2}$ be the coefficients of $C^{\gamma}$ in the expansion of $H_{(R_{\tau_1}, e_1)}(C)$ and $H_{(R_{\tau_2}, e_2)}(C)$ respectively. There is a point $(\mu_{12}, \gamma)$ on the edge $e$ which corresponds to the term of $H_{(R_{\tau}, e)}(C)$ of degree $\gamma$. We know by Lemma~\ref{momop} that 
$$R(y+\psi_i)=\sum_{0\leq j\leq n}\frac{1}{j!}R^{(j)}(\psi_i)r_j \, \, \textrm{for}\, \, i=1,2.$$
Then $ord (R^{(\gamma)}(\psi_i)-\gamma !\overline{\xi}_ix^{\mu_{12}})<\mu_{12}$ for $i=1,2$. This proves the lemma by taking $\xi_i=\gamma !\overline{\xi}_i$ for $i=1,2$. $\Box$
\newline
\newline
Let $\{\Psi_1, \dots ,\Psi_n\}$ be a fundamental system of solutions of the linear differential equation $S(y)=0$ (see e.g.~\cite{Was, Del1, Gri9}) and $\psi_1, \dots , \psi_n$ be their logarithmic derivatives respectively, i.e., $\psi_1=\Psi_1'/\Psi_1, \dots , \psi_n=\Psi_n'/\Psi_n$. Then $R(\psi_i)=0$ for all $1\leq i\leq n$. 

\dfntnl{}
Let $\psi$ be an element of the field $\mathcal{L}$. We denote by $span_r(\psi)$ the $r$-differential span of $\psi$, i.e., $span_r(\psi)$ is the $\mathbb{Z}$-module generated by $r_1(\psi), r_2(\psi), \dots $.
\edfntn

\lmml{}
Let $\psi \in \mathcal{L}$ be a solution of a {\it Riccatti} equation $R_2(y)=0$ where $R_2\in \mathbb{Z}[y_0, \dots , y_n]$ of degree $n$. Then $span_r(\psi)$ is the $\mathbb{Z}$-module generated by $r_1(\psi), \dots , r_{n-1}(\psi)$.
\elmm
\proof
Write $R_2$ in the form $R_2= r_n+ \alpha_{n-1}r_{n-1}+\cdots + \alpha_1r_1+ \alpha_0$ where $\alpha_i\in \mathbb{Z}$ for all $0\leq i< n$. Then 
\begin{eqnarray*}
r_{n+1}(\psi) & = & \psi r_n(\psi)+ Dr_n(\psi) \\
& = & \sum_{0\leq i<n}\alpha_i \big (\psi r_i(\psi)+ Dr_i(\psi)\big ) \\
& = & \sum_{0\leq i<n}\alpha_ir_{i+1}(\psi) \\
& = &   \sum_{0\leq i<n}\beta_ir_i(\psi)
\end{eqnarray*}
for suitable $\beta_i \in \mathbb{Z}$ using the fact that $r_n(\psi)=\sum_{0\leq i<n}\alpha_ir_i(\psi)$. $\Box$
\newline

Consider a $\mathbb{Z}$-module $M:=span_r(\psi_1, \dots , \psi_n)$, i.e., $M$ is the $\mathbb{Z}$-module generated by $r_1(\psi_i), r_2(\psi_i), \dots $ for all $1\leq i\leq n$. We define now what we call a $r$-cyclic vector for $M$ (this definition is similar to that of the cyclic vectors in~\cite{Kov, Chu, Van}).

\dfntnl{}
An element $m\in M$ is called a $r$-cyclic vector for $M$ if $M=span_r(m)$. 
\edfntn

The following theorem is called a $r$-cyclic vector theorem. It is similar to the cyclic vector theorem of~\cite{Kat, Kov, Chu, Van}. 

\thrml{}
Let $M$ be the $\mathbb{Z}$-module defined as above. There is a $r$-cyclic vector $m$ for $M$.
\ethrm

\crllrl{}
Let $m\in M$ be a $r$-cyclic vector for $M$. Then for any $1\leq i\leq n$, there exists a {\it Riccatti} differential polynomial $R_i\in \mathbb{Z}[y_0, \dots , y_n]$ such that $\psi_i=R_i(m)$.
\ecrllr

\lmml{ass52}
For each element $m \in M$, one can compute a {\it Riccatti} differential polynomial $R_m\in K[x][y_0, \dots , y_n]$ such that $R_m(m)=0$. In addition, there is a positive integer $s$ such that the order of $R_m(y)=0$ and the degree of $R_m$ w.r.t. $y_0, \dots , y_n$ are $\leq n^s$.  

\elmm
\proof
Each element $m \in M$ has the form $m=\alpha_1\psi_1+\cdots +\alpha_n\psi_n$ where $\alpha_1, \dots ,\alpha_n\in \mathbb{Q}$. Then 
$$m=\frac{(\Psi_1^{\alpha_1}\cdots \Psi_n^{\alpha_n})'}{\Psi_1^{\alpha_1}\cdots \Psi_n^{\alpha_n}}$$
is the logarithmic derivative of $\Psi_1^{\alpha_1}\cdots \Psi_n^{\alpha_n}$. Or Lemma 3.8 (a) of~\cite{Sin} (see also~\cite{Sin1, Van}) proves that one can construct a linear differential equation $S_m(y)=0$, denoted by $S^{\circledS^{\alpha_1+\cdots +\alpha_n}}(y)=0$ of order $\leq n^{\alpha_1+\cdots +\alpha_n}$ such that $\Psi_1^{\alpha_1}\cdots \Psi_n^{\alpha_n}$ is one of its solutions. The equation $S^{\circledS^{\alpha_1+\cdots +\alpha_n}}(y)=0$ is called the $(\alpha_1+\cdots +\alpha_n)$-th symmetric power of the linear differential equation $S(y)=0$. In order to compute the equation $S_m(y)=0$ associated with the linear combination $m=\alpha_1\psi_1+\cdots +\alpha_n\psi_n\in M$, we take the change of variable $z=y^{\alpha_1+\cdots +\alpha_n}$ where $y$ is a solution of $S(y)=0$ and we compute the successive derivatives of $z$ until we get a linear dependent family over $K$. The relation between these successive derivatives gives us the linear differential equation $S_m(z)=0$. Let $R_m$ be the {\it Riccatti} differential polynomial associated with $S_m(y)=0$, then $m$ is a solution of the equation $R_m(y)=0$. $\Box$

\rmrkl{}
For any $1\leq i\leq n$, we can take $R_{\psi_i}=R$ where $R_{\psi_i}$ is defined in Lemma~\ref{ass52}. 
\ermrk

\end{document}